\documentclass[11pt,a4paper]{article}

\usepackage{theorem,enumerate}
\usepackage{amsmath,latexsym,amssymb,amsfonts}
\usepackage{eucal}
\usepackage{color}
\usepackage{comment}

\theorembodyfont{\normalfont\slshape}
\newtheorem{thm}{Theorem}

\newtheorem{prop}[thm]{Proposition}

\newtheorem{claim}{Claim}[section]

\newtheorem{Thm}{Theorem}

\newcommand{\proof}{\medbreak\noindent\textit{Proof.}\quad}

\newcommand{\qed}{{$\quad\square$\vs{3.6}}}

\numberwithin{equation}{section}

\newcommand{\vs}[1]{\vspace*{#1 mm}}

\def\AA{{ \mathcal{A}}}
\def\BB{{ \mathcal{B}}}

\bfseries\normalfont

\title{Refinements of degree conditions for the existence\\of a spanning tree without small degree stems}

\author{
Michitaka Furuya$^{1}$\footnote{\texttt{e-mail:michitaka.furuya@gmail.com}}\and \
Akira Saito$^{2}$\footnote{\texttt{e-mail:saitou.akira@nihon-u.ac.jp}}\and \
Shoichi Tsuchiya$^{3}$\footnote{\texttt{e-mail:s.tsuchiya@isc.senshu-u.ac.jp}} \vs{5}\\
$^{1}$\textsl{College of Liberal Arts and Sciences,}\\
\textsl{Kitasato University,}\\
\textsl{1-15-1 Kitasato, Minami-ku, Sagamihara, Kanagawa 252-0373, Japan}\\
$^{2}$\textsl{Department of Information Science,}\\
\textsl{Nihon University,}\\ 
\textsl{3-25-40 Sakurajosui, Setagaya-Ku, Tokyo 156-8550, Japan}\\
$^{3}$\textsl{School of Network and Information,}\\
\textsl{Senshu University,}\\
\textsl{2-1-1 Higashimita, Tama-ku, Kawasaki-shi, Kanagawa, 214-8580, Japan}\\
}

\date{}

\begin{document}

\maketitle

\begin{abstract}
A spanning tree of a graph without no vertices of degree $2$ is called a {\it homeomorphically irreducible spanning tree} (or a {\it HIST}) of the graph.
Albertson, Berman, Hutchinson and Thomassen~[J. Graph Theory {\bf 14} (1990), 247--258] gave a minimum degree condition for the existence of a HIST, and recently, Ito and Tsuchiya~[J. Graph Theory {\bf 99} (2022), 162--170] found a sharp degree-sum condition for the existence of a HIST.
In this paper, we refine these results,
and extend the first one to a spanning tree in which no vertex other than the endvertices has small degree.
\end{abstract}



\section{Introduction}\label{sec1}

Let $G$ be a graph.
We let $V(G)$ and $E(G)$ denote the {\it vertex set} and the {\it edge set} of $G$, respectively.
For $u\in V(G)$, let $N_{G}(u)$ and $d_{G}(u)$ denote the {\it neighborhood} and the {\it degree} of $u$, respectively; thus $N_{G}(u)=\{v\in V(G):uv\in E(G)\}$ and $d_{G}(u)=|N_{G}(u)|$.
For an integer $i\geq 0$, let $V_{i}(G)=\{u\in V(G):d_{G}(u)=i\}$ and $V_{\geq i}(G)=\{u\in V(G):d_{G}(u)\geq i\}$.
We let $\delta (G)$ denote the {\it minimum degree} of $G$.
We let
$$
\sigma _{2}(G)=\min\{d_{G}(u)+d_{G}(v):u,v\in V(G),~u\neq v,~uv\notin E(G)\}
$$
if $G$ is not complete; we let $\sigma _{2}(G)=\infty $ if $G$ is complete.

Let $G$ be a connected graph.
A spanning tree of $G$ without a vertex of degree~$2$
is called a {\it homeomorphically irreducible spanning tree},
which we often abbreviate to a HIST,
of $G$.
In a sense,
a HIST is the counterpart of a hamiltonian path,
which is a spanning tree where every vertex other than the endvertices,
referring to it as a {\it stem} hereafter,
has degree~$2$.
Partly for this reason,
as the study of a hamiltonian path has evolved,
a HIST has also been attracting attention
(for example, see \cite{ABHT,CRS,CS,FT2,HNO,ZWHY}).
Also,
a HIST is worth studying in its own right.
For example,
a HIST appears as an essential part in the construction of a Halin graph,
which was constructed by Halin~\cite{H}, and was named by  Lov\'asz and Plummer~\cite{LP}.

Since the existence of a hamiltonian path is often studied
in relation with degree conditions,
it is natural to study the relationship between the existence
of a HIST and degrees of graphs.
The first significant discovery in this direction
was made by Albertson et al.~\cite{ABHT}.
An immediate consequence of Dirac's Theorem states that
a graph of order $n$ with minimum degree at least $\frac{1}{2}(n-1)$ 
contains a hamiltonian path.
While this result is sharp,
they proved that a much weaker minimum degree condition
guarantees the existence of a HIST.

\begin{Thm}[Albertson, Berman, Hutchinson and Thomassen~\cite{ABHT}]\label{ThmA}
Let $G$ be a connected graph of order $n$, and suppose that $\delta (G)\geq 4\sqrt{2n}$.
Then $G$ has a HIST.
\end{Thm}

Dirac's Theorem was later extended by Ore.
Again an immediate consequence of Ore's Theorem
says that every graph $G$ of order $n$ with $\sigma_2(G)\ge n-1$
contains a hamiltonian path.
As for a HIST,
Ito and the third author of this paper~\cite{IT}
recently proved the following result.

\begin{Thm}[Ito and Tsuchiya~\cite{IT}]\label{ThmB}
Let $G$ be a connected graph of order $n\geq 8$.
If $\sigma _{2}(G)\geq n-1$, then $G$ has a HIST.
\end{Thm}

One peculiar aspect is that
they also proved the sharpness of the result~:~They
constructed infinitely many connected graphs $G$
with $\sigma_2(G)=|V(G)|-2$
which does not contain a HIST.
Therefore,
while there is a wide discrepancy between minimum degree conditions
for the existence of a hamiltonian path and a HIST,
the gap disappears when we consider the degree sum conditions.

As a hamiltonian path has become a popular topic,
its extended notion has also been pursued.
As we have mentioned,
a hamiltonian path is a spanning tree in which every stem has degree~$2$.
Relaxing this restraint,
for a integer $k\ge 2$,
we define a {\it spanning $k$-tree}
as a spanning tree in which every stem has degree lying between $2$ and $k$.
A number of sufficient conditions for the existence
of a hamiltonian path have been naturally generalized to those for the existence
of a spanning $k$-tree.

If a HIST is the counterpart of a hamiltonian path
and a spanning $k$-tree is a natural extension of a hamiltonian path,
we may seek for a similar extension of a HIST.
Partly under this motivation,
the first and the third authors of this paper
defined $[2, k]$-ST in~\cite{FT}.
In a spanning $k$-tree,
we allow stems to have degree ranging in the interval~$[2, k]$.
So as a possible extension of a HIST,
they  forbade the degrees of stems to fall into~$[2, k]$~:~A $[2,k]$-ST
of a connected graph $G$ is a spanning tree of $G$
in which every stem has degree larger than~$k$.
Note that a $[2,1]$-ST is a spanning tree without restriction 
\footnote{For a connected graph of order $n$, we regard a spanning $(n-1)$-tree as a spanning tree without restriction.}.

Though the term $[2, k]$-ST was introduced in~\cite{FT},
its notion had implicitly appeared in the literature before it.
A graph of order $n$ is said to be {\it pancyclic}
if it contains a cycle of order $k$ for every integer $k$ with $3\le k\le n$.
As we have remarked,
a Halin graph is constructed from a HIST.
Though it is not always pancyclic,
Bondy and Lov\'asz~\cite{BL} and, independently Skowro\'nska~\cite{Sko} proved that
a Halin graph constructed from a $[2,3]$-ST is pancyclic.
Also,
when Albertson et al.~proved Theorem~\ref{ThmA},
they remarked that by adopting the same proof technique,
they could prove that for every integer~$k\ge 2$,
there exists a constant $c_k$ such that
every connected graph of order $n$ and minimum degree at least $c_k\sqrt{n}$
contains a $[2,k]$-ST.
Since the existence of a $[2,k]$-ST becomes stronger as $k$ grows,
the constant $c_k$ is at least a non-decreasing function of $k$,
and we expect that it is actually an increasing function.
However,
they did not mention the actual behavior of $c_k$.

With these backgrounds,
we set two goals in this paper.
First,
we study the behavior of the constant $c_k$ as a function of $k$
and prove $c_k=O\bigl(k^{\frac{3}{2}}\bigr)$.

\begin{thm}
\label{thm1}
Let $k \ge 1$ be an integer.
Let $G$ be a connected graph of order $n$, and suppose that $\delta (G)\geq \sqrt{k(k-1)(k+2\sqrt{2k}+2)n}$.
Then $G$ has a $[2,k]$-ST.
\end{thm}

If $k=2$, then $\sqrt{k(k-1)(k+2\sqrt{2k}+2)}=4$.
Thus Theorem~\ref{thm1} slightly improves Theorem~\ref{ThmA}.
In Section~\ref{sec-pf-thm1}, we also show $c_{k}=\Omega(k^{\frac{1}{2}})$.

The second goal is to refine Theorem~\ref{thm2}.
As we have seen,
while the minimum degree condition for the existence of a HIST
is considerably weaker than that of a hamiltonian path,
we cannot observe such a gap for the degree sum conditions.
As the first step to understand this seemingly peculiar difference,
we characterize the graphs
with $\sigma_2(G)=|V(G)|-2$
that does not contain a HIST.
Let $D_{n}$ be the graph obtained from a complete graph $K$ of order $n-2$ and a path $P$ of order $3$ by identifying a vertex of $K$ and an endvertex of $P$.
Then, as mentioned in \cite{IT}, $|V(D_{n})|=n$, $\sigma _{2}(D_{n})=n-2$ and $D_{n}$ has no HIST.

\begin{thm}
\label{thm2}
Let $G$ be a connected graph of order $n\geq 10$, and suppose that $\sigma _{2}(G)\geq n-2$.
Then $G$ has a HIST if and only if $G$ is not isomorphic to $D_{n}$.
\end{thm}

We remark that in~\cite{FT3},
the first and the third authors of this paper
extended Theorem~\ref{thm2}
to $[2,k]$-STs
by adopting a different approach.
However,
the proof strategy in~\cite{FT3}
works only for graphs of sufficiently
large order.
In particular,
the application of the result to the ordinary HIST
verifies Theorem~\ref{thm2}
only for graphs of order at least~$295$.
Many arguments concerning HISTs are made by
induction on the order of the graph,
and the number~$295$ is rather too large
to use as an induction basis.
For this reason,
we give Theorem~\ref{thm2}
in this paper.


\subsection{Further notations and a preliminary}\label{sec-pre}

In this subsection, we prepare some notations and introduce a known result, which are used in our proof.
For terms and symbols not defined in this paper, we refer the reader to \cite{D}.

Let $G$ be a graph.
For two disjoint subsets $U_{1}$ and $U_{2}$ of $V(G)$, let $E_{G}(U_{1},U_{2})=\{u_{1}u_{2}\in E(G):u_{1}\in U_{1},~u_{2}\in U_{2}\}$.
For $u,v\in V(G)$, the {\it distance} between $u$ and $v$, denoted by ${\rm dist}_{G}(u,v)$, is the minimum length of a path of $G$ connecting $u$ and $v$.
The value ${\rm diam}(G)=\max\{{\rm dist}_{G}(u,v):u,v\in V(G)\}$ is called the {\it diameter} of $G$.
For $F\subseteq E(G)$, let $V(F)=\{u,v:uv\in F\}$.
For a subgraph $H$ of $G$ and a subset $F$ of $E(G)$, let $H+F$ be the subgraph of $G$ with $V(H+F)=V(H)\cup V(F)$ and $E(H+F)=E(H)\cup F$.
Let $\omega (G)$ be the number of components of $G$.


Let $m\geq 1$ be an integer.
Let $p_{1},p_{2},\ldots ,p_{m}$ be integers with $p_{i}\geq 1$.
For each integer $i$ with $1\leq i\leq m$, let $A_{i}$ be a copy of the complete bipartite graph $K_{2,p_{i}}$, and let $x_{i,1}$ and $x_{i,2}$ be two vertices of $A_{i}$ such that $\{x_{i,1},x_{i,2}\}$ is one of the partite sets of $A_{i}$.
Let $\hat{A}_{m}(p_{1},p_{2},\ldots ,p_{m})$ be the graph obtained from $A_{1},A_{2},\ldots ,A_{m}$ by identifying $x_{1,1},x_{2,1},\ldots ,x_{m,1}$ and adding the edge set $\{x_{i,2}x_{j,2}:1\leq i<j\leq m\}$.
For an integer $p\geq 1$, let $B_{p}$ be the graph obtained from $\hat{A}_{2}(2,p)$ by adding the edge $yy'$ where $\{y,y'\}=V(A_{1})\setminus \{x_{1,1},x_{1,2}\}$.
Let $\AA =\{\hat{A}_{m}(p_{1},p_{2},\ldots ,p_{m}):m\geq 1,~p_{i}\geq 1,~1\leq i\leq m\}$ and $\BB =\{B_{p}:p\geq 1\}$ (see Figure~\ref{f-AB}).
Recently, Shan and Tsuchiya~\cite{ST} proved the following theorem.

\begin{figure}
\begin{center}
\input{fig-AB.tex}
\caption{Graphs belong to $\AA \cup \BB $}
\label{f-AB}
\end{center}
\end{figure}

\begin{Thm}[Shan and Tsuchiya~\cite{ST}]
\label{ThmC}
Let $G$ be a graph of order $n\geq 10$ with ${\rm diam}(G)=2$.
Then $G$ has a HIST if and only if $G$ is not isomorphic to any graph in $\AA \cup \BB$.
\end{Thm}

\section{Proof of Theorem~\ref{thm1}}\label{sec-pf-thm1}

\medbreak\noindent\textit{Proof of Theorem~\ref{thm1}.}\quad

It is trivial that Theorem~\ref{thm1} also holds for $k=1$. 
Thus we may assume that $k \ge 2$.
Throughout this proof, we implicitly use the fact that $k\geq 2$.
Let $c_{k}=\sqrt{k(k-1)(k+2\sqrt{2k}+2)}$.
Hence
\begin{align*}
c_{k} &= \sqrt{k(k-1)(k+2\sqrt{2k}+2)}\\
&\geq \sqrt{k(k-1)(k+4+2)}\\
&= \sqrt{k^{3}+5k^{2}-6k}.
\end{align*}
Since $n>\delta (G)\geq c_{k}\sqrt{n}$, we have $\sqrt{n}(\sqrt{n}-c_{k})>0$.
Consequently, we obtain
\begin{align}
\sqrt{n}>c_{k}\geq \sqrt{k^{3}+5k^{2}-6k}.\label{cond-1-1}
\end{align}
Let $p=\frac{(c_{k}^{2}-k(k-1)(k-2))\sqrt{n}+kc_{k}}{2c_{k}}$.

\begin{claim}
\label{cl-1-1}
We have $c_{k}\sqrt{n}>p>k+2$.
\end{claim}
\proof
By (\ref{cond-1-1}),
\begin{align*}
c_{k}\sqrt{n}-p &= c_{k}\sqrt{n}-\frac{(c_{k}^{2}-k(k-1)(k-2))\sqrt{n}+kc_{k}}{2c_{k}}\\
&= \frac{(c_{k}^{2}+k(k-1)(k-2))\sqrt{n}-kc_{k}}{2c_{k}}\\
&> \frac{(k^{3}+5k^{2}-6k+k(k-1)(k-2))c_{k}-kc_{k}}{2c_{k}}\\
&= \frac{(2k^{3}+2k^{2}-4k)-k}{2}\\
&> 0
\end{align*}
and
\begin{align*}
p-(k+2) &= \frac{(c_{k}^{2}-k(k-1)(k-2))\sqrt{n}+kc_{k}}{2c_{k}}-(k+2)\\
&> \frac{(k^{3}+5k^{2}-6k-k(k-1)(k-2))c_{k}+kc_{k}-2c_{k}(k+2)}{2c_{k}}\\
&= \frac{(8k^{2}-8k)+k-2(k+2)}{2}\\
&> 0.
\end{align*}
Therefore, $c_{k}\sqrt{n}>p>k+2$.
\qed

Let $u_{0}\in V(G)$, and let $F_{0}$ be the graph with $V(F_{0})=\{u_{0}\}\cup N_{G}(u_{0})$ and $E(F_{0})=\{u_{0}v:v\in N_{G}(u_{0})\}$.
Then by Claim~\ref{cl-1-1}, $|V(F_{0})|=d_{G}(u_{0})+1>\delta (G)\geq c_{k}\sqrt{n}>p$ and $d_{F_{0}}(u_{0})=d_{G}(u_{0})\geq c_{k}\sqrt{n}>k+2$.
In particular, there exists a subforest $F$ of $G$ such that $\bigcup _{2\leq i\leq k+1}V_{i}(F)=\emptyset $ and every component of $F$ has at least $p$ vertices.
Choose $F$ so that $|V(F)|$ is as large as possible.
Note that
\begin{align}
\omega (F)\leq \frac{n}{p}.\label{cond-1-2}
\end{align}

Let $A=V(G)\setminus V(F)$.
Fix a vertex $a\in A$.
It follows from Claim~\ref{cl-1-1} that if $|N_{G}(a)\cap A|\geq p$, then the graph $F'$ with $V(F')=V(F)\cup \{a\}\cup (N_{G}(a)\cap A)$ and $E(F')=E(F)\cup \{ay:y\in N_{G}(a)\cap A\}$ is a subforest of $G$ such that $\bigcup _{2\leq i\leq k+1}V_{i}(F')=\emptyset $ and every component of $F'$ has at least $p$ vertices, which contradicts the maximality of $F$.
Thus $|N_{G}(a)\cap A|<p$.
Again by the maximality of $F$, we also obtain $N_{G}(a)\cap V_{\geq k+2}(F)=\emptyset $.
Consequently, we have
\begin{align}
|N_{G}(a)\cap V_{1}(F)|=|N_{G}(a)\cap V(F)|=d_{G}(a)-|N_{G}(a)\cap A|>c_{k}\sqrt{n}-p.\label{cond-1-3}
\end{align}
Furthermore, for a vertex $u\in V_{1}(F)$, it follows from the maximality of $F$ that $|N_{G}(u)\cap A|\leq k$.
This together with (\ref{cond-1-3}) implies that
\begin{align*}
(c_{k}\sqrt{n}-p)|A| &< \sum _{a\in A}|N_{G}(a)\cap V_{1}(F)|\\
&= |E_{G}(A,V_{1}(F))|\\
&= \sum _{u\in V_{1}(F)}|N_{G}(u)\cap A|\\
&\leq k|V_{1}(F)|\\
&< k(n-|A|),
\end{align*}
and hence
\begin{align}
|A|<\frac{kn}{c_{k}\sqrt{n}-p+k}.\label{cond-1-4}
\end{align}

For each vertex $a\in A$, take $u_{a}\in N_{G}(a)\cap V_{1}(F)$ (here, the condition that $N_{G}(a)\cap V_{1}(F)\neq \emptyset $ is assured by (\ref{cond-1-3}) and Claim~\ref{cl-1-1}).
Let $F^{*}$ be the spanning subgraph of $G$ with $E(F^{*})=E(F)\cup \{au_{a}:a\in A\}$.
Then $F^{*}$ is a forest with $\omega (F^{*})=\omega (F)$.
Furthermore, since $\bigcup _{2\leq i\leq k+1}V_{i}(F^{*})\subseteq \{u_{a}:a\in A\}$, we have $|\bigcup _{2\leq i\leq k+1}V_{i}(F^{*})|\leq |A|$.
Since $G$ is connected and $F^{*}$ is a spanning subforest of $G$, there exists a spanning tree $T_{0}$ of $G$ with $E(F^{*})\subseteq E(T_{0})$.
Since $|E(T_{0})\setminus E(F^{*})|$ is less than $\omega (F^{*})~(=\omega (F))$, $|(\bigcup _{2\leq i\leq k+1}V_{i}(T_{0}))\setminus (\bigcup _{2\leq i\leq k+1}V_{i}(F^{*}))|\leq 2|E(T_{0})\setminus E(F^{*})|<2\omega (F)$.
In particular, it follows from (\ref{cond-1-2}) and (\ref{cond-1-4}) that
$$
\left|\bigcup _{2\leq i\leq k+1}V_{i}(T_{0})\right|<|A|+2\omega (F)<\frac{kn}{c_{k}\sqrt{n}-p+k}+\frac{2n}{p}=\frac{kpn+2(c_{k}\sqrt{n}-p+k)n}{(c_{k}\sqrt{n}-p+k)p}.
$$

For a spanning tree $S$ of $G$, let
$$
\mu (S)=\sum _{2\leq i\leq k+1}(k+2-i)|V_{i}(S)|.
$$
A spanning tree $S$ is {\it admissible} if $\mu (S)<\frac{k(kpn+2(c_{k}\sqrt{n}-p+k)n)}{(c_{k}\sqrt{n}-p+k)p}$.
Since $\mu (T_{0})\leq k\sum _{2\leq i\leq k+1}|V_{i}(T_{0})|<\frac{k(kpn+2(c_{k}\sqrt{n}-p+k)n)}{(c_{k}\sqrt{n}-p+k)p}$, $T_{0}$ is admissible.
Now we choose an admissible spanning tree $T$ of $G$ so that 
\begin{enumerate}[{\bf (X1)}]
\item[{\bf (X1)}]
$(|V_{2}(T)|,|V_{3}(T)|,\ldots ,|V_{k}(T)|)$ is lexicographically as small as possible.
\end{enumerate}
We show that $T$ is a $[2,k]$-ST of $G$.
By way of contradiction, suppose that $\bigcup _{2\leq i\leq k}V_{i}(T)\neq \emptyset $.
Let $s$ be the minimum integer such that $2\leq s\leq k$ and $V_{s}(T)\neq \emptyset $, and let $u\in V_{s}(T)$.

\begin{claim}
\label{cl-1-2}
We have $d_{G}(u)>(k-1)\mu (T)$.
\end{claim}
\proof
Since $d_{G}(u)\geq c_{k}\sqrt{n}$ and $T$ is an admissible spanning tree of $G$, it suffices to show that
\begin{align}
c_{k}\sqrt{n}\geq \frac{(k-1)k(kpn+2(c_{k}\sqrt{n}-p+k)n)}{(c_{k}\sqrt{n}-p+k)p}.\label{target-cl}
\end{align}

Since $x=c_{k}^{2}~(=k(k-1)(k+2\sqrt{2k}+2))$ is a solution to an equation $x^{2}-2k(k-1)(k+2)x+k^{2}(k-1)^{2}(k-2)^{2}=0$, we have
$$
c_{k}^{4}-2k(k-1)(k+2)c_{k}^{2}+k^{2}(k-1)^{2}(k-2)^{2}=0.
$$
Furthermore, it follows from (\ref{cond-1-1}) that
$$
c_{k}^{2}-k^{3}-k^{2}+2k\geq (k^{3}+5k^{2}-6k)-k^{3}-k^{2}+2k>0.
$$
Hence
\begin{align*}
&((c_{k}^{2}-k(k-1)(k-2))\sqrt{n}+kc_{k})^{2}-8c_{k}k(k-1)(c_{k}\sqrt{n}+k)\sqrt{n}\\
&= (c_{k}^{4}-2k(k-1)(k+2)c_{k}^{2}+k^{2}(k-1)^{2}(k-2)^{2})n+2kc_{k}(c_{k}^{2}-k^{3}-k^{2}+2k)\sqrt{n}+k^{2}c_{k}^{2}\\
&> 0,
\end{align*}
and so
\begin{align}
\frac{((c_{k}^{2}-k(k-1)(k-2))\sqrt{n}+kc_{k})^{2}}{4c_{k}}-2k(k-1)(c_{k}\sqrt{n}+k)\sqrt{n}>0.\label{cond-cl-2-2-3}
\end{align}
Since $p=\frac{(c_{k}^{2}-k(k-1)(k-2))\sqrt{n}+kc_{k}}{2c_{k}}$, 
\begin{align*}
&-c_{k}\sqrt{n}p^{2}+((c_{k}^{2}-k(k-1)(k-2))n+kc_{k}\sqrt{n})p\\
&= -c_{k}\sqrt{n}\left(p-\frac{(c_{k}^{2}-k(k-1)(k-2))\sqrt{n}+kc_{k}}{2c_{k}}\right)^{2}+\frac{((c_{k}^{2}-k(k-1)(k-2))\sqrt{n}+kc_{k})^{2}\sqrt{n}}{4c_{k}}\\
&= \frac{((c_{k}^{2}-k(k-1)(k-2))\sqrt{n}+kc_{k})^{2}\sqrt{n}}{4c_{k}}.
\end{align*}
This together with (\ref{cond-cl-2-2-3}) leads to
\begin{align*}
c_{k}\sqrt{n}(c_{k}\sqrt{n}&-p+k)p-(k-1)k(kpn+2(c_{k}\sqrt{n}-p+k)n)\\
&= -c_{k}\sqrt{n}p^{2}+((c_{k}^{2}-k(k-1)(k-2))n+kc_{k}\sqrt{n})p-2k(k-1)(c_{k}\sqrt{n}+k)n\\
&= \left(\frac{((c_{k}^{2}-k(k-1)(k-2))\sqrt{n}+kc_{k})^{2}}{4c_{k}}-2k(k-1)(c_{k}\sqrt{n}+k)\sqrt{n}\right)\sqrt{n}\\
&> 0.
\end{align*}
Consequently, (\ref{target-cl}) holds.
\qed

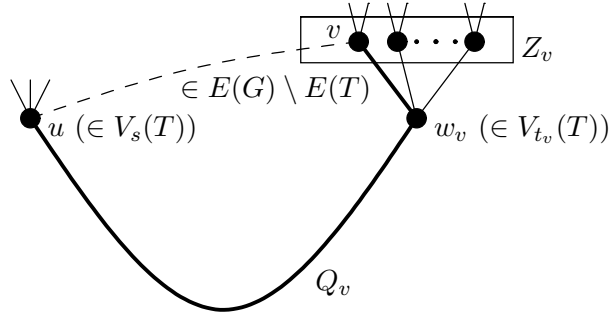
\begin{figure}
\begin{center}
{\unitlength 0.1in%
\begin{picture}(26.4000,16.0000)(11.0000,-20.0000)%
%
\special{sh 1.000}%
\special{ia 1200 1000 50 50 0.0000000 6.2831853}%
\special{pn 8}%
\special{ar 1200 1000 50 50 0.0000000 6.2831853}%
%
\special{sh 1.000}%
\special{ia 3200 1000 50 50 0.0000000 6.2831853}%
\special{pn 8}%
\special{ar 3200 1000 50 50 0.0000000 6.2831853}%
%
\special{sh 1.000}%
\special{ia 2900 600 50 50 0.0000000 6.2831853}%
\special{pn 8}%
\special{ar 2900 600 50 50 0.0000000 6.2831853}%
%
\special{sh 1.000}%
\special{ia 3100 600 50 50 0.0000000 6.2831853}%
\special{pn 8}%
\special{ar 3100 600 50 50 0.0000000 6.2831853}%
%
\special{sh 1.000}%
\special{ia 3500 600 50 50 0.0000000 6.2831853}%
\special{pn 8}%
\special{ar 3500 600 50 50 0.0000000 6.2831853}%
%
\special{pn 8}%
\special{pa 1100 800}%
\special{pa 1200 1000}%
\special{fp}%
\special{pa 1200 800}%
\special{pa 1200 1000}%
\special{fp}%
\special{pa 1200 1000}%
\special{pa 1300 800}%
\special{fp}%
%
\special{pn 8}%
\special{pa 2900 600}%
\special{pa 2850 400}%
\special{fp}%
\special{pa 2950 400}%
\special{pa 2950 400}%
\special{fp}%
%
\special{pn 8}%
\special{pa 2950 400}%
\special{pa 2900 600}%
\special{fp}%
%
\special{pn 8}%
\special{pa 3100 600}%
\special{pa 3050 400}%
\special{fp}%
\special{pa 3150 400}%
\special{pa 3150 400}%
\special{fp}%
%
\special{pn 8}%
\special{pa 3150 400}%
\special{pa 3100 600}%
\special{fp}%
%
\special{pn 8}%
\special{pa 3500 600}%
\special{pa 3450 400}%
\special{fp}%
\special{pa 3550 400}%
\special{pa 3550 400}%
\special{fp}%
%
\special{pn 8}%
\special{pa 3550 400}%
\special{pa 3500 600}%
\special{fp}%
%
\special{pn 20}%
\special{pa 3200 1000}%
\special{pa 2900 600}%
\special{fp}%
%
\special{pn 8}%
\special{pa 3100 600}%
\special{pa 3200 1000}%
\special{fp}%
%
\special{pn 8}%
\special{pa 3200 1000}%
\special{pa 3500 600}%
\special{fp}%
\put(27.7000,-18.6000){\makebox(0,0){$Q_{v}$}}%
\put(12.9000,-9.8000){\makebox(0,0)[lt]{$u~(\in V_{s}(T))$}}%
%
\special{pn 8}%
\special{pa 2600 470}%
\special{pa 3700 470}%
\special{pa 3700 710}%
\special{pa 2600 710}%
\special{pa 2600 470}%
\special{pa 3700 470}%
\special{fp}%
\put(27.7000,-5.5000){\makebox(0,0){$v$}}%
\put(37.4000,-7.0000){\makebox(0,0)[lb]{$Z_{v}$}}%
\put(32.9000,-9.8000){\makebox(0,0)[lt]{$w_{v}~(\in V_{t_{v}}(T))$}}%
\put(24.7000,-8.5000){\makebox(0,0){$\in E(G)\setminus E(T)$}}%
%
\special{pn 8}%
\special{pn 8}%
\special{pa 1200 1000}%
\special{pa 1230 989}%
\special{pa 1261 979}%
\special{fp}%
\special{pa 1323 957}%
\special{pa 1352 947}%
\special{pa 1383 937}%
\special{fp}%
\special{pa 1445 915}%
\special{pa 1503 896}%
\special{pa 1506 895}%
\special{fp}%
\special{pa 1569 875}%
\special{pa 1595 866}%
\special{pa 1625 856}%
\special{pa 1630 855}%
\special{fp}%
\special{pa 1692 835}%
\special{pa 1748 819}%
\special{pa 1754 817}%
\special{fp}%
\special{pa 1817 799}%
\special{pa 1871 783}%
\special{pa 1879 781}%
\special{fp}%
\special{pa 1942 764}%
\special{pa 1994 751}%
\special{pa 2004 749}%
\special{fp}%
\special{pa 2068 733}%
\special{pa 2087 729}%
\special{pa 2119 722}%
\special{pa 2131 720}%
\special{fp}%
\special{pa 2195 706}%
\special{pa 2212 703}%
\special{pa 2244 697}%
\special{pa 2258 694}%
\special{fp}%
\special{pa 2323 683}%
\special{pa 2338 680}%
\special{pa 2370 675}%
\special{pa 2386 672}%
\special{fp}%
\special{pa 2451 661}%
\special{pa 2464 659}%
\special{pa 2496 655}%
\special{pa 2514 652}%
\special{fp}%
\special{pa 2579 643}%
\special{pa 2591 641}%
\special{pa 2623 636}%
\special{pa 2643 634}%
\special{fp}%
\special{pa 2708 625}%
\special{pa 2718 623}%
\special{pa 2771 616}%
\special{fp}%
\special{pa 2836 608}%
\special{pa 2846 607}%
\special{pa 2877 603}%
\special{pa 2900 600}%
\special{fp}%
%
\special{pn 4}%
\special{sh 1}%
\special{ar 3200 600 10 10 0 6.2831853}%
\special{sh 1}%
\special{ar 3400 600 10 10 0 6.2831853}%
\special{sh 1}%
\special{ar 3300 600 10 10 0 6.2831853}%
\special{sh 1}%
\special{ar 3300 600 10 10 0 6.2831853}%
%
\special{pn 20}%
\special{pa 3200 1000}%
\special{pa 3177 1034}%
\special{pa 3155 1068}%
\special{pa 3132 1102}%
\special{pa 3109 1135}%
\special{pa 3087 1169}%
\special{pa 3064 1202}%
\special{pa 3042 1236}%
\special{pa 3019 1269}%
\special{pa 2996 1301}%
\special{pa 2974 1334}%
\special{pa 2951 1366}%
\special{pa 2928 1397}%
\special{pa 2906 1429}%
\special{pa 2883 1459}%
\special{pa 2861 1490}%
\special{pa 2838 1519}%
\special{pa 2815 1549}%
\special{pa 2793 1577}%
\special{pa 2770 1605}%
\special{pa 2747 1632}%
\special{pa 2725 1659}%
\special{pa 2702 1685}%
\special{pa 2680 1710}%
\special{pa 2657 1735}%
\special{pa 2634 1758}%
\special{pa 2612 1781}%
\special{pa 2566 1823}%
\special{pa 2544 1843}%
\special{pa 2521 1862}%
\special{pa 2499 1880}%
\special{pa 2453 1912}%
\special{pa 2431 1926}%
\special{pa 2408 1940}%
\special{pa 2385 1952}%
\special{pa 2363 1962}%
\special{pa 2340 1972}%
\special{pa 2318 1980}%
\special{pa 2295 1987}%
\special{pa 2272 1992}%
\special{pa 2250 1996}%
\special{pa 2227 1999}%
\special{pa 2204 2000}%
\special{pa 2182 2000}%
\special{pa 2159 1998}%
\special{pa 2137 1994}%
\special{pa 2114 1989}%
\special{pa 2091 1983}%
\special{pa 2069 1975}%
\special{pa 2046 1966}%
\special{pa 2023 1956}%
\special{pa 2001 1944}%
\special{pa 1978 1932}%
\special{pa 1955 1918}%
\special{pa 1933 1902}%
\special{pa 1910 1886}%
\special{pa 1888 1869}%
\special{pa 1842 1831}%
\special{pa 1820 1811}%
\special{pa 1774 1767}%
\special{pa 1752 1744}%
\special{pa 1729 1720}%
\special{pa 1707 1695}%
\special{pa 1661 1643}%
\special{pa 1639 1616}%
\special{pa 1593 1560}%
\special{pa 1571 1531}%
\special{pa 1548 1501}%
\special{pa 1526 1471}%
\special{pa 1503 1441}%
\special{pa 1480 1409}%
\special{pa 1458 1378}%
\special{pa 1412 1314}%
\special{pa 1390 1281}%
\special{pa 1367 1248}%
\special{pa 1345 1215}%
\special{pa 1322 1182}%
\special{pa 1299 1148}%
\special{pa 1277 1115}%
\special{pa 1231 1047}%
\special{pa 1209 1013}%
\special{pa 1200 1000}%
\special{fp}%
\end{picture}}%
\caption{A vertex $v$ and related notations}
\label{f1}
\end{center}
\end{figure}

Fix a vertex $v\in N_{G}(u)\setminus N_{T}(u)$.
Let $Q_{v}$ be the unique path in $T$ joining $u$ and $v$, and write $N_{Q_{v}}(v)=\{w_{v}\}$.
Let $t_{v}$ be the integer with $w_{v}\in V_{t_{v}}(T)$.
Since $uv\notin E(T)$, $w_{v}\neq u$.
In particular, $|N_{Q_{v}}(w_{v})|=2$.
By the definition of $s$, we have
\begin{align}
2\leq s\leq t_{v}.\label{cond-1-tv}
\end{align}
Let $Z_{v}=\{v\}\cup (N_{T}(w_{v})\setminus N_{Q_{v}}(w_{v}))$ (see Figure~\ref{f1}).
Then $|Z_{v}|=t_{v}-1$.
Remark that
\begin{align}
\mbox{for vertices $v,v'\in N_{G}(u)\setminus N_{T}(u)$, $v'\in Z_{v}$ if and only if $w_{v'}=w_{v}$.}\label{cond-1-Zv}
\end{align}

\begin{claim}
\label{cl-1-3}
Let $v\in N_{G}(u)\setminus N_{T}(u)$.
Then
\begin{enumerate}[{\upshape(i)}]
\item
$t_{v}\in \{s,s+1\}$, and
\item
$Z_{v}\not \subseteq N_{G}(u)$.
\end{enumerate}
\end{claim}
\proof
We first prove (i).
Suppose that $t_{v}\notin \{s,s+1\}$.
By (\ref{cond-1-tv}), $t_{v}\geq s+2$.
Let $T_{1}=(T-vw_{v})+uv$.
Then $T_{1}$ is a spanning tree of $G$, and $d_{T_{1}}(u)=d_{T}(u)+1=s+1$, $d_{T_{1}}(w_{v})=d_{T}(w_{v})-1=t_{v}-1~(\geq s+1)$ and $d_{T_{1}}(y)=d_{T}(y)$ for all vertices $y\in V(G)\setminus \{u,w_{v}\}$.
Now we calculate the value $\mu (T)-\mu (T_{1})$.
Recall that $s\leq k$.
If $t_{v}\leq k+1$, then
$$
\mu (T)-\mu (T_{1}) = (k+2-s)+(k+2-t_{v})-(k+2-(s+1))-(k+2-(t_{v}-1))=0;
$$
if $t_{v}=k+2$, then
$$
\mu (T)-\mu (T_{1}) = (k+2-s)-(k+2-(s+1))-(k+2-(k+1))=0;
$$
if $t_{u}\geq k+3$, then
$$
\mu (T)-\mu (T_{1}) = (k+2-s)-(k+2-(s+1))=1.
$$
In any case, we have $\mu (T)\geq \mu (T_{1})$, and hence $T_{1}$ is an admissible spanning tree of $G$.
However, by the definition of $s$, $|V_{i}(T_{1})|=|V_{i}(T)|=0$ for all integers $i$ with $2\leq i\leq s-1$ and $|V_{s}(T_{1})|=|V_{s}(T)\setminus \{u\}|=|V_{s}(T)|-1$, which contradicts (X1).
Thus (i) holds.

Next we prove (ii).
Suppose that $Z_{v}\subseteq N_{G}(u)$.
Let $T_{2}=(T-\{w_{v}z:z\in Z_{v}\})+\{uz:z\in Z_{v}\}$.
Then $T_{2}$ is a spanning tree of $G$, and $d_{T_{2}}(u)=d_{T}(u)+(t_{v}-1)=s+t_{v}-1~(\geq s+1)$, $d_{T_{2}}(w_{v})=1$ and $d_{T_{2}}(y)=d_{T}(y)$ for all vertices $y\in V(G)\setminus \{u,w_{v}\}$.
Now we calculate the value $\mu (T)-\mu (T_{2})$.
By (i), $t_{v}\leq s+1\leq k+1$.
If $s+t_{v}-1\leq k+1$, then
$$
\mu (T)-\mu (T_{2}) = (k+2-s)+(k+2-t_{v})-(k+2-(s+t_{v}-1))=k+1;
$$
if $s+t_{v}-1\geq k+2$, then
$$
\mu (T)-\mu (T_{2}) = (k+2-s)+(k+2-t_{v})>0.
$$
In either case, we have $\mu (T)>\mu (T_{2})$, and hence $T_{2}$ is an admissible spanning tree of $G$.
However, by the definition of $s$, $|V_{i}(T_{2})|=|V_{i}(T)|=0$ for all integers $i$ with $2\leq i\leq s-1$ and $|V_{s}(T_{2})|=|V_{s}(T)\setminus \{u\}|=|V_{s}(T)|-1$, which contradicts (X1).
Thus (ii) holds.
\qed

By Claim~\ref{cl-1-3}(i), $N_{G}(u)\setminus N_{T}(u)=\{v\in N_{G}(u)\setminus N_{T}(u):w_{v}\in V_{s}(T)\cup V_{s+1}(T)\}$.
By (\ref{cond-1-Zv}) and Claim~\ref{cl-1-3}(ii),
$$
\mbox{for a vertex $w\in V_{s}(T)$, $|\{v\in N_{G}(u)\setminus N_{T}(u): w_{v}=w\}|\leq s-2$}
$$
and
$$
\mbox{for a vertex $w\in V_{s+1}(T)$, $|\{v\in N_{G}(u)\setminus N_{T}(u): w_{v}=w\}|\leq s-1$.}
$$
Again by (\ref{cond-1-Zv}), it follows that
\begin{align}
|N_{G}(u)\setminus N_{T}(u)| &= \sum _{w\in V_{s}(T)\cup V_{s+1}(T)}|\{v\in N_{G}(u)\setminus N_{T}(u): w_{v}=w\}|\nonumber \\
&\leq (s-2)|V_{s}(T)|+(s-1)|V_{s+1}(T)|.\label{cond-1-7}
\end{align}
Since
$$
(k-1)(k+2-s)-(2s-2) = k^{2}+k-s(k+1)\geq k^{2}+k-k(k+1)=0
$$
and
$$
(k-1)(k+2-(s+1))-(s-1) = k^{2}-sk\geq k^{2}-k^{2}=0,
$$
we have
\begin{align}
\mbox{$(k-1)(k+2-s)\geq 2s-2$~~and~~$(k-1)(k+2-(s+1))\geq s-1$.}\label{cond-1-8}
\end{align}
Since $u\in V_{s}(T)$, $|V_{s}(T)|\geq 1$.
This together with Claim~\ref{cl-1-2}, (\ref{cond-1-7}) and (\ref{cond-1-8}) implies that
\begin{align*}
d_{G}(u) &> (k-1)\mu (T)\\
&\geq (k-1)((k+2-s)|V_{s}(T)|+(k+2-(s+1))|V_{s+1}(T)|)\\
&\geq (2s-2)|V_{s}(T)|+(s-1)|V_{s+1}(T)|\\
&= s|V_{s}(T)|+(s-2)|V_{s}(T)|+(s-1)|V_{s+1}(T)|\\
&\geq s+|N_{G}(v)\setminus N_{T}(u)|,
\end{align*}
and hence $d_{G}(u)-s>|N_{G}(u)\setminus N_{T}(u)|$, which contradicts the fact that $u\in V_{s}(T)$.

This completes the proof of Theorem~\ref{thm1}.
\qed

We expect that the coefficient of $\sqrt{n}$ in Theorem~\ref{thm1} can be further improved.
However, in our proof, the coefficient of $\sqrt{n}$ is best possible:
Let $c$ be a positive number that is (slightly) smaller than $\sqrt{k(k-1)(k+2\sqrt{2k}+2)}$, and suppose that a connected graph $G$ of order $n$ satisfies $\delta (G)=c\sqrt{n}$.
Then we can verify that $c^{4}-2k(k-1)(k+2)c^{2}+k^{2}(k-1)^{2}(k-2)^{2}<0$, and hence (\ref{cond-cl-2-2-3}) does not hold if $n$ is sufficiently large.

\begin{figure}
\begin{center}
{\unitlength 0.1in%
\begin{picture}(24.0000,20.0000)(3.0000,-29.0000)%
%
{\color[named]{Black}{%
\special{pn 20}%
\special{ar 600 1200 300 300 0.0000000 6.2831853}%
}}%
%
{\color[named]{Black}{%
\special{pn 20}%
\special{ar 1400 1200 300 300 0.0000000 6.2831853}%
}}%
%
{\color[named]{Black}{%
\special{pn 20}%
\special{ar 2400 1200 300 300 0.0000000 6.2831853}%
}}%
%
{\color[named]{Black}{%
\special{pn 0}%
\special{sh 1.000}%
\special{ia 600 1400 40 40 0.0000000 6.2831853}%
}}%
{\color[named]{Black}{%
\special{pn 8}%
\special{ar 600 1400 40 40 0.0000000 6.2831853}%
}}%
%
{\color[named]{Black}{%
\special{pn 0}%
\special{sh 1.000}%
\special{ia 1400 1400 40 40 0.0000000 6.2831853}%
}}%
{\color[named]{Black}{%
\special{pn 8}%
\special{ar 1400 1400 40 40 0.0000000 6.2831853}%
}}%
%
{\color[named]{Black}{%
\special{pn 0}%
\special{sh 1.000}%
\special{ia 2400 1400 40 40 0.0000000 6.2831853}%
}}%
{\color[named]{Black}{%
\special{pn 8}%
\special{ar 2400 1400 40 40 0.0000000 6.2831853}%
}}%
%
{\color[named]{Black}{%
\special{pn 0}%
\special{sh 1.000}%
\special{ia 2400 2000 40 40 0.0000000 6.2831853}%
}}%
{\color[named]{Black}{%
\special{pn 8}%
\special{ar 2400 2000 40 40 0.0000000 6.2831853}%
}}%
%
{\color[named]{Black}{%
\special{pn 0}%
\special{sh 1.000}%
\special{ia 1400 2000 40 40 0.0000000 6.2831853}%
}}%
{\color[named]{Black}{%
\special{pn 8}%
\special{ar 1400 2000 40 40 0.0000000 6.2831853}%
}}%
%
{\color[named]{Black}{%
\special{pn 0}%
\special{sh 1.000}%
\special{ia 600 2000 40 40 0.0000000 6.2831853}%
}}%
{\color[named]{Black}{%
\special{pn 8}%
\special{ar 600 2000 40 40 0.0000000 6.2831853}%
}}%
%
{\color[named]{Black}{%
\special{pn 20}%
\special{pa 300 1900}%
\special{pa 2700 1900}%
\special{pa 2700 2100}%
\special{pa 300 2100}%
\special{pa 300 1900}%
\special{pa 2700 1900}%
\special{fp}%
}}%
%
{\color[named]{Black}{%
\special{pn 20}%
\special{pa 300 2300}%
\special{pa 2700 2300}%
\special{pa 2700 2500}%
\special{pa 300 2500}%
\special{pa 300 2300}%
\special{pa 2700 2300}%
\special{fp}%
}}%
%
{\color[named]{Black}{%
\special{pn 20}%
\special{pa 300 2700}%
\special{pa 2700 2700}%
\special{pa 2700 2900}%
\special{pa 300 2900}%
\special{pa 300 2700}%
\special{pa 2700 2700}%
\special{fp}%
}}%
\put(15.0000,-22.0000){\makebox(0,0){{\color[named]{Black}{$+$}}}}%
\put(15.0000,-26.0000){\makebox(0,0){{\color[named]{Black}{$+$}}}}%
%
{\color[named]{Black}{%
\special{pn 8}%
\special{pa 600 2000}%
\special{pa 600 1400}%
\special{fp}%
}}%
%
{\color[named]{Black}{%
\special{pn 8}%
\special{pa 1400 1400}%
\special{pa 1400 2000}%
\special{fp}%
}}%
%
{\color[named]{Black}{%
\special{pn 8}%
\special{pa 2400 2000}%
\special{pa 2400 1400}%
\special{fp}%
}}%
%
{\color[named]{Black}{%
\special{pn 4}%
\special{sh 1}%
\special{ar 1800 1600 16 16 0 6.2831853}%
\special{sh 1}%
\special{ar 2000 1600 16 16 0 6.2831853}%
\special{sh 1}%
\special{ar 1900 1600 16 16 0 6.2831853}%
\special{sh 1}%
\special{ar 1900 1600 16 16 0 6.2831853}%
}}%
\put(6.0000,-11.0000){\makebox(0,0){{\color[named]{Black}{$H_{1}$}}}}%
\put(14.0000,-11.0000){\makebox(0,0){{\color[named]{Black}{$H_{2}$}}}}%
\put(24.0000,-11.0000){\makebox(0,0){{\color[named]{Black}{$H_{l}$}}}}%
\put(7.5000,-20.0000){\makebox(0,0){{\color[named]{Black}{$u_{1}$}}}}%
\put(15.5000,-20.0000){\makebox(0,0){{\color[named]{Black}{$u_{2}$}}}}%
\put(25.5000,-20.0000){\makebox(0,0){{\color[named]{Black}{$u_{l}$}}}}%
\put(28.5000,-24.0000){\makebox(0,0){{\color[named]{Black}{$U_{2}$}}}}%
\put(34.2000,-28.0000){\makebox(0,0){{\color[named]{Black}{$U_{1}\setminus \{u_{i}:1\leq i\leq l\}$}}}}%
\put(7.3000,-13.4000){\makebox(0,0){{\color[named]{Black}{$v_{1}$}}}}%
\put(15.3000,-13.4000){\makebox(0,0){{\color[named]{Black}{$v_{2}$}}}}%
\put(25.3000,-13.4000){\makebox(0,0){{\color[named]{Black}{$v_{l}$}}}}%
\end{picture}}%
\caption{Graph $G_{k,d}$}
\label{f-Gkd}
\end{center}
\end{figure}
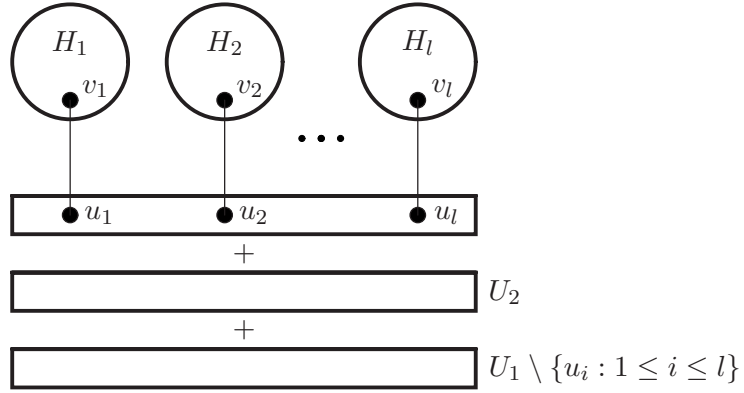

We conclude this section by showing that the coefficient of $\sqrt{n}$ in Theorem~\ref{thm1} cannot be improved with $\sqrt{k-1}-\varepsilon $ for any $\varepsilon >0$.
It suffices to construct infinitely many connected graphs $G$ of order $n$ such that $\delta (G)\geq \sqrt{(k-1)n+\frac{(2k-1)^{2}}{4}}-k+\frac{1}{2}$ and $G$ has no $[2,k]$-ST.
Let $k$ and $d$ be positive integers such that $k\geq 2$ and $l:=\frac{d}{k-1}$ is an integer.
Let $H_{0},H_{1},\ldots ,H_{l}$ be vertex-disjoint graphs such that $H_{0}$ is a complete bipartite graph having a bipartition $(U_{1},U_{2})$ with $|U_{1}|=|U_{2}|=d$ and for an integer $i$ with $1\leq i\leq l$, $H_{i}$ is a complete graph of order $d+1$.
Take $l$ vertices $u_{1},u_{2},\ldots ,u_{l}\in U_{1}$, and for each integer $i$ with $1\leq i\leq l$, take a vertex $v_{i}\in V(H_{i})$.
Let $G_{k,d}=(\bigcup _{0\leq i\leq l}H_{i})+\{u_{i}v_{i}:1\leq i\leq l\}$ (see Figure~\ref{f-Gkd}).
Then $G_{k,d}$ is a connected graph of order $n:=2d+\frac{d(d+1)}{k-1}$ and $\delta (G_{k,d})=d=\frac{\sqrt{4(k-1)n+(2k-1)^{2}}-2k+1}{2}$.
This together with the following proposition leads to a desired conclusion.

\begin{prop}
\label{prop-Gkd}
There exists no $[2,k]$-ST of $G_{k,d}$.
\end{prop}
\proof
Suppose that $G_{k,d}$ has a $[2,k]$-ST $T$.
Let $G_{0}$ be the subgraph of $T$ induced by $\{u_{i}:1\leq i\leq l\}\cup U_{2}$.
Then $|V(G_{0})|=l+d=\frac{kd}{k-1}$.
For each integer $i$ with $1\leq i\leq l$, since $u_{i}$ is a cut-vertex of $G_{k,d}$, $d_{T}(u_{i})\geq k+1$.
Hence $|E(G_{0})|=\sum _{1\leq i\leq l}|N_{T}(u_{i})\cap U_{2}|\geq kl=\frac{kd}{k-1}=|V(G_{0})|$.
This implies that $G_{0}$ contains a cycle, which contradicts the assumption that $G_{0}$ is a subgraph of a tree $T$.
\qed

\section{Proof of Theorem~\ref{thm2}}\label{sec-pf-thm2}

Throughout this proof, we implicitly use the fact that $n\geq 10$.
As we mentioned in Section~\ref{sec1}, $D_{n}$ has no HIST, and so the ``only if'' part holds.
Thus we suppose that that $G$ is not isomorphic to $D_{n}$ and prove that $G$ has a HIST.

Set $d=\mbox{diam}(G)$.
If $d=1$, then $G$ is a complete graph, and hence $G$ has a HIST.
Since we can easily verify that every graph $A\in \AA\cup \BB$ of order $n$ has two non-adjacent vertices of degree $2$, we have $\sigma _{2}(A)=4~(<n-2)$.
In particular, $G$ is isomorphic to no graph in $\AA\cup \BB$.
Hence by Theorem~\ref{ThmC}, if $d=2$, then $G$ has a HIST.
Thus we may assume that $d\geq 3$.

Let $P$ be a diametral path of $G$, i.e., $P$ is a path of length $d$ joining two vertices with the distance $d$ in $G$, and write $P=u_{0}u_{1}\cdots u_{d}$.
Choose $P$ so that
\begin{enumerate}[{\bf (P1)}]
\setlength{\parskip}{0cm}
\setlength{\itemsep}{0cm}
\item
$d_{G}(u_{0})$ is as small as possible,
\item
subject to (P1), $\max\{3-d_{G}(u_{1}),0\}$ is as small as possible, and
\item
subject to (P2), $|N_{G}(u_{0})\cap (N_{G}(u_{1})\cup N_{G}(u_{2}))|$ is as large as possible.
\end{enumerate}
Since $d\geq 3$, $N_{G}(u_{0})$ and $N_{G}(u_{d})$ are disjoint, and hence
$$
n-2\leq \sigma _{2}(G)\leq d_{G}(u_{0})+d_{G}(u_{d})=|N_{G}(u_{0})\cup N_{G}(u_{d})|\leq |V(G)\setminus \{u_{0},u_{d}\}|=n-2.
$$
This forces $N_{G}(u_{0})\cup N_{G}(u_{d})=V(G)\setminus \{u_{0},u_{d}\}$ and $\sigma _{2}(G)=d_{G}(u_{0})+d_{G}(u_{d})=n-2$.
Since $u_{2}\notin N_{G}(u_{0})$, we have $u_{2}\in N_{G}(u_{d})$, and so $d=3$.
In the remaining of the proof, we implicitly use the facts that $N_{G}(u_{0})\cup N_{G}(u_{3})=V(G)\setminus \{u_{0},u_{3}\}$ and $\sigma _{2}(G)=d_{G}(u_{0})+d_{G}(u_{3})=n-2$.

\begin{claim}
\label{cl-2-1}
Let $i\in \{0,3\}$, and let $v\in N_{G}(u_{i})$.
Then
\begin{enumerate}[{\upshape(i)}]
\item
$d_{G}(v)\geq n-2-d_{G}(u_{3-i})=d_{G}(u_{i})$, and
\item
if $N_{G}(v)\cap N_{G}(u_{3-i})=\emptyset $, then $N_{G}(v)=(N_{G}(u_{i})\setminus \{v\})\cup \{u_{i}\}$, and in particular, $N_{P}(u_{i})\neq \{v\}$ and $v$ is adjacent to the unique vertex in $N_{P}(u_{i})$.
\end{enumerate}
\end{claim}
\proof
Since $vu_{3-i}\notin E(G)$, $d_{G}(v)\geq \sigma _{2}(G)-d_{G}(u_{3-i})=n-2-d_{G}(u_{3-i})=d_{G}(u_{i})$, which proves (i).

We suppose that $N_{G}(v)\cap N_{G}(u_{3-i})=\emptyset $, and prove (ii).
Since $N_{G}(v)\subseteq (N_{G}(u_{i})\setminus \{v\})\cup \{u_{i}\}$, we have $d_{G}(v)\leq |(N_{G}(u_{i})\setminus \{v\})\cup \{u_{i}\}|=d_{G}(u_{i})$.
This together with (i) forces $N_{G}(v)=(N_{G}(u_{i})\setminus \{v\})\cup \{u_{i}\}$.
In particular, either $N_{P}(u_{i})=\{v\}$ or $v$ is adjacent to the unique vertex in $N_{P}(u_{i})$.
If $N_{P}(u_{i})=\{v\}$, then $v$ is a vertex of $P$, and so $N_{G}(v)\cap N_{G}(u_{3-i})\neq \emptyset $, which contradicts the assumption in this paragraph.
Consequently, we obtain the desired conclusion.
\qed

Now we focus on the following conditions for a subtree $T$ of $G$:
\begin{enumerate}[{\bf (T1)}]
\setlength{\parskip}{0cm}
\setlength{\itemsep}{0cm}
\item
$\{u_{0},u_{3}\}\cup N_{G}(u_{0})\subseteq V(T)$,
\item
$V(T)\setminus \{u_{3}\}\subseteq V_{1}(T)\cup V_{\geq 3}(T)$, and
\item
$n-|V(T)|+d_{T}(u_{3})-3\geq 0$.
\end{enumerate}

\begin{claim}
\label{cl-2-2}
If there exists a subtree $T$ of $G$ satisfying (T1)--(T3), then $G$ has a HIST.
\end{claim}
\proof
Suppose that there exists a subtree $T$ of $G$ satisfying (T1)--(T3).
Since every vertex in $V(G)\setminus V(T)$ is adjacent to $u_{3}$ by (T1), $T':=T+\{u_{3}v:v\in V(G)\setminus V(T)\}$ is a spanning tree of $G$.
Furthermore, it follows from (T3) that $|V(G)\setminus V(T)|=n-|V(T)|\geq 3-d_{T}(u_{3})$, and so 
$$
d_{T'}(u_{3}) = d_{T}(u_{3})+|V(G)\setminus V(T)|\geq d_{T}(u_{3})+(3-d_{T}(u_{3}))=3.
$$
This together with (T2) implies that $T'$ is a HIST of $G$.
\qed

We divide the proof into two cases.

\medskip
\noindent
\textbf{Case 1:} $d_{G}(u_{0})\leq 3$.

\begin{claim}
\label{cl-2-3}
We have $d_{G}(u_{1})\geq 3$.
\end{claim}
\proof
Suppose that $d_{G}(u_{1})\leq 2$.
Then $N_{G}(u_{1})=\{u_{0},u_{2}\}$.
For the moment, we further suppose that $d_{G}(u_{0})=1$.
For a vertex $v\in N_{G}(u_{3})\setminus \{u_{2}\}$, since $N_{G}(v)\cap N_{G}(u_{0})=\emptyset $, it follows from Claim~\ref{cl-2-1}(ii) with $i=3$ that $N_{G}(v)=(N_{G}(u_{3})\setminus \{v\})\cup \{u_{3}\}$.
Since $v$ is arbitrary, $N_{G}(u_{3})\cup \{u_{3}\}$ is a clique of $G$.
This implies that $G$ is isomorphic to $D_{n}$, which is a contradiction.
Thus $d_{G}(u_{0})\geq 2$.

Let $v'\in N_{G}(u_{0})\setminus \{u_{1}\}$.
Since $d_{G}(u_{1})=2$, $u_{1}v'\notin E(G)$.
Hence by Claim~\ref{cl-2-1}(ii) with $(i,v)=(0,v')$, we have $N_{G}(v')\cap N_{G}(u_{3})\neq \emptyset $.
In particular, $P'=u_{0}v'wu_{3}$ is a diametral path of $G$ where $w\in N_{G}(v')\cap N_{G}(u_{3})$.
Since $u_{1}v'\notin E(G)$, $d_{G}(v')+2=d_{G}(v')+d_{G}(u_{1})\geq \sigma _{2}(G)\geq n-2>4$, and hence $d_{G}(v')\geq 3$.
Since $P'$ is a diametral path of $G$, this contradicts to (P2).
\qed

By Claim~\ref{cl-2-3}, there exists a vertex $a\in N_{G}(u_{1})\setminus \{u_{0},u_{2}\}$.
Choose $a$ so that $a\in N_{G}(u_{0})$ if possible.

\begin{claim}
\label{cl-2-4}
If $N_{G}(u_{0})\subseteq N_{G}(u_{1})\cup N_{G}(u_{2})$, then $G$ has a HIST.
\end{claim}
\proof
Suppose that $N_{G}(u_{0})\subseteq N_{G}(u_{1})\cup N_{G}(u_{2})$.
For each $v\in N_{G}(u_{0})\setminus \{u_{1}\}$, take $w_{v}\in N_{G}(v)\cap \{u_{1},u_{2}\}$.
By Claim~\ref{cl-2-1}(i) with $(i,v)=(3,u_{2})$, $d_{G}(u_{2})\geq n-2-d_{G}(u_{0})\geq 10-2-3=5$, and hence there exists a vertex $a'\in N_{G}(u_{2})\setminus \{u_{1},u_{3},a\}$.
Let $T=P+(\{u_{1}a,u_{2}a'\}\cup \{vw_{v}:v\in N_{G}(u_{0})\setminus \{u_{1},a,a'\}\})$.
Then $T$ is a tree and satisfies (T1) and (T2).
Furthermore, $n-|V(T)|+d_{T}(u_{3})-3\geq n-8+1-3\geq 0$, and so $T$ satisfies (T3).
Hence by Claim~\ref{cl-2-2}, $G$ has a HIST.
\qed

By Claim~\ref{cl-2-4}, we may assume that $N_{G}(u_{0})\not \subseteq N_{G}(u_{1})\cup N_{G}(u_{2})$.
Since $u_{1}\in N_{G}(u_{2})$, this implies that $d_{G}(u_{0})\geq 2$ and there exists a vertex $z_{1}\in N_{G}(u_{0})\setminus (N_{G}(u_{1})\cup N_{G}(u_{2}))$.
In particular,
\begin{align}
|N_{G}(u_{0})\cap (N_{G}(u_{1})\cup N_{G}(u_{2}))|\leq d_{G}(u_{0})-1.\label{cond-2-1}
\end{align}
Since $u_{1}z_{1}\notin E(G)$, it follows from Claim~\ref{cl-2-1}(ii) with $(i,v)=(0,z_{1})$ that $N_{G}(z_{1})\cap N_{G}(u_{3})\neq \emptyset $.
Let $z_{2}\in N_{G}(z_{1})\cap N_{G}(u_{3})$.
Since $u_{2}z_{1}\notin E(G)$, $z_{2}\neq u_{2}$.

\begin{claim}
\label{cl-2-5}
If $N_{G}(u_{0})\setminus \{u_{1},z_{1}\}\subseteq N_{G}(u_{1})$, then $u_{1}z_{2}\notin E(G)$, and in particular, $z_{2}\neq a$.
\end{claim}
\proof
Suppose that $N_{G}(u_{0})\setminus \{u_{1},z_{1}\}\subseteq N_{G}(u_{1})$ and $u_{1}z_{2}\in E(G)$.
Then $u_{0}u_{1}z_{2}u_{3}$ is a diametral path of $G$.
However, it follows from (\ref{cond-2-1}) that $|N_{G}(u_{0})\cap (N_{G}(u_{1})\cup N_{G}(z_{2}))|=|N_{G}(u_{0})|=d_{G}(u_{0})>|N_{G}(u_{0})\cap (N_{G}(u_{1})\cup N_{G}(u_{2}))|$, which contradicts (P3).
\qed

By Claim~\ref{cl-2-1}(i) with $(i,v)=(3,u_{2})$, $d_{G}(u_{2})\geq n-2-d_{G}(u_{0})\geq 5$.
Since $u_{2}z_{1}\notin E(G)$, this implies that there exist vertices $b,b'\in N_{G}(u_{2})$ such that $b\notin \{u_{1},z_{1},a,z_{2},u_{3}\}$ and $b'\notin N_{G}(u_{0})\cup \{a,u_{3}\}$ where $b$ might be equal to $b'$.

\begin{figure}
\begin{center}
\input{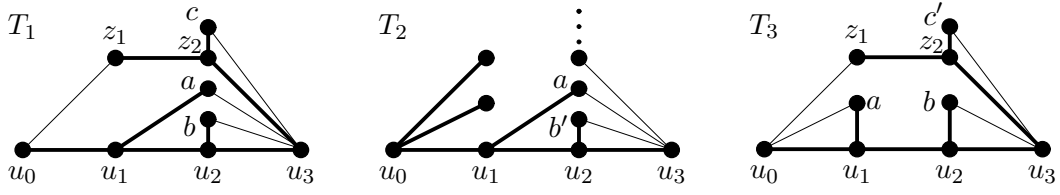}
\caption{Trees $T_{1}$, $T_{2}$ and $T_{3}$}
\label{f-T1}
\end{center}
\end{figure}

Suppose that $d_{G}(u_{0})=2$.
Then by Claim~\ref{cl-2-1}(i) with $(i,v)=(3,z_{2})$, $d_{G}(z_{2})\geq n-2-d_{G}(u_{0})\geq 6$, and hence there exists a vertex $c\in N_{G}(z_{2})\setminus \{z_{1},a,u_{2},b,u_{3}\}$.
By Claim~\ref{cl-2-5}, $z_{2}\neq a$ and $c\neq u_{1}$.
Hence by the definition of $b$ and $c$, the vertices $u_{0},u_{1},u_{2},u_{3},z_{1},z_{2},a,b,c$ are pairwise distinct.
Let $T_{1}=P+\{u_{1}a,u_{2}b,u_{3}z_{2},z_{2}z_{1},z_{2}c\}$ (see the left graph in Figure~\ref{f-T1}).
Then $T_{1}$ is a tree and satisfies (T1) and (T2).
Furthermore, $n-|V(T_{1})|+d_{T_{1}}(u_{3})-3\geq n-9+2-3\geq 0$, and so $T_{1}$ satisfies (T3).
Hence by Claim~\ref{cl-2-2}, $G$ has a HIST.
Thus we may assume that $d_{G}(u_{0})=3$.

Suppose that $a\notin N_{G}(u_{0})$.
Then by the choice of $a$, $N_{G}(u_{0})\cap N_{G}(u_{1})=\emptyset $.
Recall that $a\neq b'$ and $b'\notin N_{G}(u_{0})$.
Let $T_{2}=P+(\{u_{1}a,u_{2}b'\}\cup \{u_{0}v:v\in N_{G}(u_{0})\setminus \{u_{1}\}\})$ (see the middle graph in Figure~\ref{f-T1}).
Then $T_{2}$ is a tree and satisfies (T1) and (T2).
Furthermore, $n-|V(T_{2})|+d_{T_{2}}(u_{3})-3\geq n-8+1-3\geq 0$, and so $T_{2}$ satisfies (T3).
Hence by Claim~\ref{cl-2-2}, $G$ has a HIST.
Thus we may assume that $a\in N_{G}(u_{0})$, i.e., $N_{G}(u_{0})=\{u_{1},z_{1},a\}$.

If $z_{2}a\in E(G)$, then $u_{0}az_{2}u_{3}$ is a diametral path of $G$, $\max\{3-d_{G}(a),0\}=0$ and $|N_{G}(u_{0})\cap (N_{G}(a)\cup N_{G}(z_{2}))|=3>|N_{G}(u_{0})\cap (N_{G}(u_{1})\cup N_{G}(u_{2}))|$ by (\ref{cond-2-1}), which contradicts (P3).
Thus
\begin{align}
z_{2}a\notin E(G).\label{cond-2-2}
\end{align}
By Claim~\ref{cl-2-1}(i) with $(i,v)=(3,z_{2})$, $d_{G}(z_{2})\geq n-2-d_{G}(u_{0})\geq 5$, and hence there exists a vertex $c'\in N_{G}(z_{2})\setminus \{z_{1},u_{2},b,u_{3}\}$.
By Claim~\ref{cl-2-5} and (\ref{cond-2-2}), $c'\notin \{u_{1},a\}$.
Hence by the definition of $b$ and $c'$, the vertices $u_{0},u_{1},u_{2},u_{3},z_{1},z_{2},a,b,c'$ are pairwise distinct.
Let $T_{3}=P+\{u_{1}a,u_{2}b,u_{3}z_{2},z_{2}z_{1},z_{2}c'\}$ (see the right graph in Figure~\ref{f-T1}).
Then $T_{3}$ is a tree and satisfies (T1) and (T2).
Furthermore, $n-|V(T_{3})|+d_{T_{3}}(u_{3})-3\geq n-9+2-3\geq 0$, and so $T_{3}$ satisfies (T3).
Hence by Claim~\ref{cl-2-2}, $G$ has a HIST.

\medskip
\noindent
\textbf{Case 2:} $d_{G}(u_{0})\geq 4$.

By (P1), $d_{G}(u_{3})\geq d_{G}(u_{0})\geq 4$.
This together with Claim~\ref{cl-2-1}(i) implies that
\begin{align}
\mbox{$d_{G}(v)\geq d_{G}(u_{i})\geq 4$ for $v\in V(G)\setminus \{u_{0},u_{3}\}$,}\label{cond-2-2+}
\end{align}
where $i$ is the integer with $i\in \{0,3\}$ and $v\in N_{G}(u_{i})$.
In particular, we can take vertices $a_{j}\in N_{G}(u_{j})\setminus V(P)~(j\in \{1,2\})$ with $a_{1}\neq a_{2}$ so that
\begin{enumerate}[{\bf (A1)}]
\setlength{\parskip}{0cm}
\setlength{\itemsep}{0cm}
\item[{\bf (A1)}]
$\min\{|N_{G}(u_{i})\setminus \{a_{1},a_{2}\}|:i\in \{0,3\}\}$ is as large as possible, and
\item[{\bf (A2)}]
subject to (A1), $|\{u_{1}a_{2},u_{2}a_{1}\}\cap E(G)|$ is as small as possible.
\end{enumerate}

If $\min\{|N_{G}(u_{i})\setminus \{a_{1},a_{2}\}|:i\in \{0,3\}\}\geq 3$, then $P+(\{u_{1}a_{1},u_{2}a_{2}\}\cup \{u_{i}v:v\in N_{G}(u_{i})\setminus \{a_{1},a_{2}\},~i\in \{0,3\}\})$ is a HIST of $G$, as desired.
Thus we may assume that there exists an integer $i_{0}\in \{0,3\}$ such that $|N_{G}(u_{i_{0}})\setminus \{a_{1},a_{2}\}|\leq 2$.
Since $d_{G}(u_{3})\geq d_{G}(u_{0})\geq 4$, it follows from (A1) that $d_{G}(u_{i_{0}})=4$, $\{a_{1},a_{2}\}\subseteq N_{G}(u_{i_{0}})\setminus V(P)$ and
\begin{align}
\mbox{$(N_{G}(u_{j})\setminus V(P))\cap N_{G}(u_{3-i_{0}})=\emptyset $ for each $j\in \{1,2\}$.}\label{cond-2-3}
\end{align}
In particular, $|N_{G}(u_{i_{0}})\cap (N_{G}(u_{1})\cup N_{G}(u_{2}))|\geq |N_{P}(u_{i_{0}})\cup \{a_{1},a_{2}\}|=3>1=|N_{P}(u_{3-i_{0}})|=|N_{G}(u_{3-i_{0}})\cap (N_{G}(u_{1})\cup N_{G}(u_{2}))|$.
Since $d_{G}(u_{3})\geq d_{G}(u_{0})$, this together with (\ref{cond-2-2+}) and the choice of $P$ (i.e., (P1)--(P3)) leads to $i_{0}=0$.

\begin{figure}
\begin{center}
\input{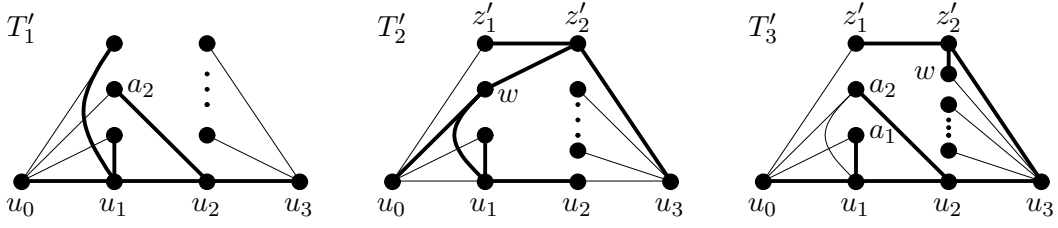}
\caption{Trees $T'_{1}$, $T'_{2}$ and $T'_{3}$}
\label{f-T2}
\end{center}
\end{figure}

If $N_{G}(u_{0})\setminus \{u_{1},a_{2}\}\subseteq N_{G}(u_{1})$, then $T'_{1}:=P+(\{u_{2}a_{2}\}\cup \{u_{1}v:v\in N_{G}(u_{0})\setminus \{u_{1},a_{2}\}\})$ is a tree satisfying (T1)--(T3) because $n-|V(T'_{1})|+d_{T'_{1}}(u_{3})-3=n-7+1-3>0$ (see the left graph in Figure~\ref{f-T2}).
Hence by Claim~\ref{cl-2-2}, $G$ has a HIST.
Thus we may assume that $N_{G}(u_{0})\setminus \{u_{1},a_{2}\}\not \subseteq N_{G}(u_{1})$.
Since $d_{G}(u_{1})\geq d_{G}(u_{0})=4$ by (\ref{cond-2-2+}), it follows from (\ref{cond-2-3}) that $|N_{G}(u_{0})\setminus (N_{G}(u_{1})\cup \{u_{1}\})|=1$, and so $N_{G}(u_{0})\cap N_{G}(u_{1})=\{a_{1},a_{2}\}$.
Write $N_{G}(u_{0})\setminus (N_{G}(u_{1})\cup \{u_{1}\})=\{z'_{1}\}$.
If $z'_{1}u_{2}\in E(G)$, then $|N_{G}(u_{0})\setminus \{a_{1},z'_{1}\}|=2=|N_{G}(u_{0})\setminus \{a_{1},a_{2}\}|$ and $|\{u_{1}z'_{1},u_{2}a_{1}\}\cap E(G)|=|\{u_{2}a_{1}\}\cap E(G)|<|\{u_{1}a_{2},u_{2}a_{1}\}\cap E(G)|$, which contradicts (A1) and (A2).
Thus $z'_{1}u_{2}\notin E(G)$.
This together with Claim~\ref{cl-2-1}(ii) with $(i,v)=(0,z'_{1})$ implies that $z'_{1}$ is adjacent to a vertex $z'_{2}\in N_{G}(u_{3})\setminus \{u_{2}\}$.
Since $d_{G}(z'_{2})\geq 4$ by (\ref{cond-2-2+}), there exists a vertex $w\in N_{G}(z'_{2})\setminus \{z'_{1},u_{3}\}$.
By (\ref{cond-2-3}), $w\notin \{u_{1},u_{2}\}$.

Suppose that $w\in N_{G}(u_{0})$.
Then $w\in \{a_{1},a_{2}\}$.
Let $T'_{2}$ be a subgraph of $G$ with $V(T'_{2})=V(G)\setminus (N_{G}(u_{3})\setminus \{u_{2},z'_{2}\})$ and $E(T'_{2})=\{u_{1}u_{2},u_{1}a_{1},u_{1}a_{2},wu_{0},wz'_{2},z'_{2}z'_{1},z'_{2}u_{3}\}$ (see the middle graph in Figure~\ref{f-T2}).
Then $T'_{2}$ is a tree and satisfies (T1) and (T2).
Furthermore, $n-|V(T'_{2})|+d_{T'_{2}}(u_{3})-3=n-8+1-3\geq 0$, and so $T'_{2}$ satisfies (T3).
Hence by Claim~\ref{cl-2-2}, $G$ has a HIST.
Thus we may assume that $w\in N_{G}(u_{3})$.
Recall that $w\neq u_{2}$.
Let $T'_{3}=P+\{u_{1}a_{1},u_{2}a_{2},u_{3}z'_{2},z'_{2}z'_{1},z'_{2}w\}$ (see the right graph in Figure~\ref{f-T2}).
Then $T'_{3}$ is a tree and satisfies (T1) and (T2).
Furthermore, $n-|V(T'_{3})|+d_{T'_{3}}(u_{3})-3=n-9+2-3\geq 0$, and so $T'_{3}$ satisfies (T3).
Hence by Claim~\ref{cl-2-2}, $G$ has a HIST.

This completes the proof of Theorem~\ref{thm2}.

\section*{Acknowledgment}

This work was supported by the Research Institute for Mathematical Sciences, an International Joint Usage/Research Center located in Kyoto University.
This work is also supported by JSPS KAKENHI Grant numbers 18K13449 (to M.F.), 23K03204 (to M.F.), 20K11684 (to A.S.) and 19K14584 (to S.T.) and research grant of Senshu University 2023 (to S.T).

\end{document}